\input amstex
\documentstyle{amsppt}
\magnification=\magstep 1
\TagsOnRight
\topmatter
\title Link Homotopy in $S^n \times \Bbb R^{m-n}$ 
and Higher Order $\mu$--Invariants \endtitle
\author Ulrich Koschorke *\endauthor
\leftheadtext{Ulrich Koschorke}
\affil Universit\"at Siegen, Emmy-Noether-Campus, D-57068 Siegen \endaffil
\address Universit\"at Siegen,
Emmy Noether Campus, Walter-Flex-Str. 3,
D-57068 Siegen, Germany
\endaddress
\email koschorke\@mathematik.uni-siegen.de \endemail
\thanks * Supported within the German-Brazilian Cooperation by IB-BMBF and CNPq.
\endthanks
\abstract
Given a suitable link map \ $f$ \ into a manifold \ $M$, \ we constructed, in \cite{10}, link homotopy invariants \ $\kappa (f)$ \ and \ $\mu (f)$. In the present paper we study the case \ $M = S^n \times \Bbb R^{m - n}$ \ in detail. Here \ $\mu (f)$ turns out to be the starting term of a whole sequence \ $\mu^{(s)}(f), \ s = 0, 1, \dots,$ \ of higher \ $\mu$-invariants which together capture all the information contained in \ $\kappa (f)$. We discuss the geometric significance of these new invariants. In several instances we obtain complete classification results. A central ingredient of our approach is the homotopy theory of wedges of spheres.
\endabstract
\keywords link map, link homotopy invariant, Hopf homomorphism, linking coefficient, stable suspension, configuration space
\endkeywords
\endtopmatter

\input boxedeps.tex
\SetRokickiEPSFSpecial
\HideDisplacementBoxes

\define\incl{\operatorname{incl}}

\define\quot{\operatorname{quot}}

\define\pinch{\operatorname{pinch}}
\define\id{\operatorname{id}}

\define\km{\kappa_{\sssize M}}

\define\krm{\kappa_{\sssize{\Bbb R^m}}}
\define\mum{\mu_{\sssize M}}
\define\mur{\mu_{\sssize{\Bbb R^m}}}
\define\mug{\mu_{{\sssize M}, \gamma}}
\define\murg{\mu_{{\sssize{\Bbb R^m}}, \gamma}}
\define\pip{\pi_{|p|}}


\def\bfmu{\hbox{$\mu$}}

\document



\specialhead 1.\ \ Introduction
\endspecialhead

Given dimensions \ $p_1, \dots, p_r \ge 1$ \ and \ $m > n \ge 1$, \ we want to classify {\it link maps}
$$
f \ = \ f_1 \ \amalg \dots \amalg \ f_r \ : \ S^{p_1} \ \amalg \dots \amalg \ S^{p_r} \ \longrightarrow \ M \ := \ S^n \times \Bbb R^{m -n}
\tag 1.1
$$
(i.e.\ the spheres \ $S^{p_j}$ \ have pairwise disjoint images under the continuous maps \ $f_j \ , \ 1 \le j \le r)$ \ up to {\it link homotopy}\ (i.e.\ continuous  deformations through such link maps). Our approach centres around the homotopy class (in the standard sense)
$$
\km (f) \ := \ [\widehat f] \ \in \ [S^{p_1} \times \dots \times S^{p_r}, \widetilde C_r (M)]
\tag 1.2
$$
where the product map \ $\widehat f \ := \ f_1 \times \dots \times f_r$ \ takes values in the configuration space of ordered $r$-tuples of pairwise distinct points in \ $M$.

Unfortunately, this very natural link homotopy invariant lies in a rather unwieldy homotopy set. However, if \ $f$ \ is \ $\km$-{\it Brunnian}\ (i.e.\ if \ $\widehat f$ \ is nulhomotopic when restricted to the complement of a point in \ $S^{p_1} \times \dots \times S^{p_r})$, one can simplify (a base point preserving version of)  \ $\km (f)$ \ considerably and then extract the \lq\lq numerical\rq\rq\ link homotopy invariant
$$
\mum (f) \ \in \ \bigoplus^{(r - 2)!} \ \pi^S_{p_1 + \dots + p_r - (r -1)(m -2) - 1}
\tag 1.3
$$
which generalizes e.g.\ Milnor's \ $\mu$-invariants of classical links and, in particular, the classical linking number when \ $r = 2$. As the example of the higher dimensional Borromean link illustrates, \ $\mum (f)$ \ is in general weaker than \ $\km (f)$ \ (cf.\ \cite{10}, 5.9).

In the present paper we measure this loss of information. In Section 3 we show that \ $\mum (f)$ \ is only the starting element of a whole sequence \ $\{ \mum^{(s)} (f)\}_{s \ge 0}$ \ of higher order \ $\mu$-invariants which together characterize \ $\km (f)$ \ (cf.\ theorem 3.5). 

Thus the question arises: what is the geometric meaning of the remaining higher invariants? We prove that in most interesting cases they are the standard \ $\mur$-invariants of the {\it augmented} \ link maps \ $f^{(s)}$ \ in \ $\Bbb R^m$ \ which constist of \ $f$, \ included into \ $\Bbb R^m$, \ together with a finite number \ $s$ \ of meridians \ $\{ z_j\} \times S^{m -n -1}$ \ around \ $M \ \cong \ S^n \times \overset\circ\to B^{m - n}$ \ (cf.\ theorems 3.8 and 4.3).

\example{Example 1.4} \ Assume \ $p_1 = \dots = p_r = m - 3 \le r \ge 2$. Let the (smoothly embedded) link \ $f : \coprod^r_{j = 1} S^{p_j} \ \hookrightarrow \ M = S^n \times \Bbb R^{m - n}$ \ be homotopy Brunnian (i.e.\ every proper sub-link is link homotopically trivial). Then \ $f$ \ is also \ $\km$-Brunnian and the invariants \ $\{ \mum^{(s)} (f)\}_{s \ge 0}$ \ determine \ $f$ \ completely up to link homotopy. In particular, \ $f$ \ itself is link homotopically trivial precisely if the component maps \ $f_j : S^{p_j} \to M \sim S^n$ \ are nulhomotopic, \ $j = 1, \dots, r$, and the \ $\krm$-invariants of all augmented link maps \ $f^{(s)}$ \ vanish, \ $s = 0, 1, 2, \dots$.

This is a very special case of corollary 5.5 below (cf.\ also 3.10). \ \hfill $\blacksquare$
\endexample

The requirement that a link map \ $f$ is \ $\km$-Brunnian allows us to concentrate on linking phenomena of highest order, but it is quite restrictive. For instance we can define \ $\mu^{(s)} (f)$ \ only if all \ $\mu$-invariants of all proper sub-link maps of \ $f$ \ are defined and trivial. However, if \ $p_1, \dots, p_r \le m -3$ \ this restriction can be avoided: \ in section 4 we define a sum operation (conceivably without additive inverses) and use it to extend our invariants to  {\it all}\ link maps whether they are \ $\km$-Brunnian or not. In particular, this allows us to introduce the {\it total higher $\mu$-invariant \ $\bfmu (f)$} \ which consists of the homotopy classes \ $[f_j] \ \in \ \pi_{p_j} (S^n)$ \ of the component maps as well as of all (higher) \ $\mum$-invariants of all sub-link maps \ $f_{j_1} \amalg \dots \amalg f_{j_s}, \ 2 \le s \le r$, of $f$.

\example\nofrills{Example 1.5:} \ \ 
$p_1 = \dots = p_r = 3, \ m = 6 > n \ge 1$. \ Given \ $r \ge 1$, let \ $BLM_{3, \dots, 3}  (M_n)$ \ denote the semigroup of all $($base point preserving if \ $n = 1)$ link homotopy classes of link maps
$$
f \ = \ \coprod^r_{j = 1} f_j \qquad : \qquad \coprod^r S^3 \ @>>{ \qquad }> \ M_n \ := \ S^n \times \Bbb R^{6 -n} \ .
$$
Then for \ $n = 1$ \ the total higher \ $\mu$-invariant \ $\bfmu$ \ is injective on
$$
BLM_{3, \dots, 3} (M_n) \ \cong \ \left(\bigoplus^\infty \Bbb Z_2\right)^{r \choose 2} \ \oplus \ \left(\bigoplus^\infty  \Bbb Z\right)^{r \choose 3} \ ;
$$
for \ $n > 1,$ \ \ $\bfmu$ \ establishes the isomorphism
$$
BLM_{3, \dots, 3} (M_n) \ \cong \ \Bbb Z_2^{r \choose 2} \ \oplus \ \Bbb Z^{r \choose 3} \ \oplus \ \Bbb Z^a \qquad \qquad \ 
$$
where \ $a = {{r + 1} \choose 2}, \ r \ \text{and} \ 0$, resp., for $n = 2, \ 3$ and $\ge 4$, resp.

Moreover, the obvious inclusions \ $M_n \subset M_{n + 1}$ \ induce homomorphisms which annihilate \ $\Bbb Z^a$ \ and correspond to summation or to the identity on the remaining direct summands. \hfill $\blacksquare$
\endexample 

This will be proved in section 5 where we compare quite generally our invariants to linking coefficients of embedded links in codimensions greater than $2$. It turns out that the transition from link isotopy to link  homotopy translates into applying Freudenthal suspensions to crucial components of the linking coefficients. (In view of the close relation between the exact link homotopy sequence and James' EHP-sequence (cf.\ \cite{8}, theorem 3.1) this comes as no surprise).

Such investigations open the way to further injectivity, surjectivity, and nontriviality results for our invariants and also to criteria deciding which aspects of the linking coefficient depend only on the link homotopy class of an embedded link. 

\example\nofrills{Example 1.6:} \ \ $n = 1, \ r = 2, \ p_1, p_2 \le m - 3$. Assume that \ $\pi_{p_2} (S^{m - p_1 - 1})$ \ is stable (i.e.\ the stable suspension is bijective). Then the total higher \ $\mu$-invariant \ $\bfmu$ \ establishes an isomorphism from the semigroup consisting of all base point preserving link homotopy classes of those link maps \ $f_1 \amalg f_2 \ : \ S^{p_1} \amalg S^{p_2} \longrightarrow S^1 \times \Bbb R^{m -1}$ \ where \ $f_1$ \ is a smooth embedding, onto \ $\pi_{p_1} (S^1) \oplus \pi_{p_2} (S^1) \oplus \bigoplus^\infty \pi^S_{p_1 + p_2 - m + 1}$. \ If \ $p_1, p_2 > 1$ \ then \ $\bfmu (f)$ \ consists essentially of the linking numbers (in the universal covering space \ $\Bbb R^m$ \ of \ $S^1 \times \Bbb R^{m -1}$) of one fixed lifting \ $\widetilde f_1$ \ of  \ $f_1$ \ with all possible liftings of \ $f_2$. This is closely related to the link concordance classification (e.g.\ of augmented links of the form \ $f^{(1)})$ \ by W.\ Mio \cite{12} where similar pairwise linking numbers occur in the off-diagonal entries of his matrices. \hfill $\blacksquare$
\endexample

If we add \lq\lq parallel longitudes\rq\rq\ to an embedded link \ $f$ \  our techniques allow us also to capture certain components of the linking coefficient of \ $f$ which are not invariant under link homotopy. Details will be given elsewhere. Let us just note here that this approach leads sometimes to a full {\it isotopy} \ classification. E.\ g.\ if \ $2 \le p \le \frac23 (m -2)$ \ it follows from \cite{3}, Corollary 1.3 (compare also \cite{5}, \S\ 4) that any smooth embedding \ $f_1 : S^p \hookrightarrow S^1 \times \Bbb R^{m - 1}$ \ which is nullisotopic when included into \ $\Bbb R^m$, is already determined up to isotopy by \ $\bfmu (f_1 \amalg f^+_1)$ \ where \ $f^+_1$ \ is such a parallel longitude. The pairwise linking numbers which essentially constitute \ $\bfmu (f_1 \amalg f_1^+)$ \ play also a central role in the link {\it concordance} \ classification of such links as \ $f_1^{(1)} : S^{p_1} \amalg S^{m -2} \hookrightarrow \Bbb R^m$ \ (cf.\ \cite{12}, 3.3).

When \ $M = S^n \times \Bbb R^{m -n}$ \ the proper setting for analyzing \ $\kappa$-invariants and linking coefficients is provided by homotopy groups of wedges of spheres. In section 2 we study these in detail, not just via the \lq\lq incoming\rq\rq\ Hilton isomorphisms, but mainly via geometrically defined \lq\lq outgoing\rq\rq\ Hopf homomorphisms which give a much better understanding  of our invariants and which enjoy many compatibility properties. Here our central results are the injectivity  theorems 2.12 and 2.27.

\subhead Notations and Conventions\endsubhead \ From now on we always assume \ $r \ge 2, \  m \ge 3, \ p_1, \dots, p_r \ge 1$. $\Sigma_s$ \ and \ $E^{(\infty)}$ \ denote the permutation group in \ $s$ \ elements and (stable) suspension, resp. All spheres and their wedges are equipped with base points. Homotopy classes in such wedges are identified, via Pontryagin-Thom, with bordisms classes of framed links.

\vskip7mm

\specialhead  2. \ \ Pinching and Hopf homomorphisms
\endspecialhead

Throughout this section let \ $n \ge 1$ \ and \ $r, \ q_1, \dots \ 
q_{r -1} \ge 2$ \ be natural numbers; put
$$
|q| \ := \ q_1 \ + \dots + \ q_{r -1} \ . 
\tag 2.1
$$

We will extend the geometric Hopf homomorphism discussed e.g. in 
\S\ 3 of \cite{9} in order to study the homotopy groups 
of the wedge
$$
W \ := \ S^n  \vee \bigvee^{r -1}_{j = 1} S_j^{q_j} 
\tag 2.2
$$
of spheres of the indicated dimensions. More specifically, in view 
of later applications we will be interested in certain subgroups 
such as the {\it reduced} \ groups
$$
\pi'_k (W) \ := \ \bigcap^{r -1}_{i = 1} \ker \left(\pi_k (W) \ 
\longrightarrow \ \pi_k 
(S^n \vee \bigvee^{r -1} \Sb j = 1 \\ j \ne i\endSb S_j^{q_j})\right), 
\tag 2.3
$$
$k \in \Bbb Z$, defined by the obvious collapsing maps.

\definition{Definition 2.4}
Given an integer \ $s \ge 0$, \ a permutation $\gamma$ of the set \ 
$\{ 1,\dots, s,\dots,$\ $r + s - 2\}$ \ is called {\it s-monotone} 
\ if \ $\gamma^{-1} (1) < \dots < \ \gamma^{- 1} (s)$.

Let \ $\sum_{r, s} \subset \sum_{r + s - 2}$ \ denote the subset
of $s$-monotone permutations in the full permutation group of \ $r 
+ s - 2$ \ elements.
\enddefinition

Since the $s$-monotone permutations form a system of 
representatives of the cosets  in \ $\sum_s \backslash \sum_{r + s 
- 2}$ \ we have
$$
u_{r,s} \ := \  |\sum_{r,s}| \ = \ (r + s - 2)! / s!
\tag 2.5
$$

Given \ $s \ge 0$, \ consider also the $s$-fold pinch map
$$
\pinch_s \ : \ (S^n, *) \ \longrightarrow \ \bigvee^s_{j = 1} 
(S^n_j, *)
$$
which has degree 1 when collapsed to any wedge summand \ $S^n_j$ \ 
(and, if \ $n = 1$, traverses the circles \ $S^n_j$ \ in the order 
given by the subindex \ $j)$; these requirements determine \ 
$\pinch_s$ \ uniquely up to homotopy. For every permutation \ 
$\gamma \in \sum_{r + s - 2}$ \ we define the homomorphism
$$
H_{s, \gamma} \ : \ \pi_k (W) \ \longrightarrow \ \pi^S_{k - sn - 
|q| + r + s - 2} \quad ,
\tag 2.6
$$
$k \in \Bbb Z$, \ by \ $H_{s, \gamma} = h_\gamma \circ (\pinch_s 
\vee \id)_*$ \ where \ $h_\gamma$ \ is the Hopf invariant 
described in \S \ 3 of \cite{9}. Summing over all \ {\it 
$s$-monotone} \ permutations \ $\gamma$ \ we obtain
$$
H_s \ := \ \bigoplus_{\gamma \in \sum_{r, s}} H_{s, \gamma} \ \ : \ 
\ \pi_k (W) \longrightarrow (\pi^S_{k - sn - |q| + r + s - 2})^{u_{r, 
s}} \quad .
\tag 2.7
$$
We define the {\it total Hopf homomorphism}  \ $H$ \ by the direct 
product
$$
H \ := \prod_{s \ge 0} H_s \ \ : \ \ \pi_k (W) \longrightarrow 
\prod_{s \ge 0} (\pi^S_{k - sn - |q| + r + s - 2})^{u_{r, s}} \ \ 
.
\tag 2.8
$$

\vskip2mm
Now consider first {\it the case when} \ $n \ge 2$ \ (then the above product 
is actually a finite direct sum).

Given \ $s \ge 0$, the homotopy class of \ $\pinch_s$ \ remains 
unchanged if we permute the wedge summands in \ $\vee^s S^n_j$; 
thus using only $s$-monotone permutations \ $\gamma$ \ in 2.7 
above helps us to avoid unnecessary redundancies.

Moreover, we choose a Hilton decomposition of \ $\pi_* (W)$
(cf.\ \cite{6}) and we define 
$$
\pi''_k (W) \ := \ \bigoplus_{t \ge 0}  \pi''_k (W, t) \ \
\subset \pi'_k (W), \ \ k \in \Bbb Z \ , 
\tag 2.9
$$
where for any integer \ $t \ge 0$ \ \ $\pi''_k (W, t)$ \
denotes the direct summand corresponding to all basic Whitehead 
products which involve the inclusion \ $\iota_0 : S^n \subset W$ 
\ precisely \ $t$ \ times and each inclusion \ $\iota_j : 
S^{q_j}_j \subset W$ \ precisely once, \ $j = 1, \dots, r -1$.

Note that any such iterated Whitehead product \ $w$, basic or not, 
is a \ $\Bbb Z$--linear combination of Whitehead products of the 
form
$$
\iota_\delta \ := \ [\iota_{\delta (1)}, [\iota_{\delta (2)}, 
[\dots, [\iota_{\delta (r + t - 2)}, \iota_{r -1}]]]]
\tag 2.10
$$
where the same factors occur the same number of times, bracketed 
in the indicated special fashion and arranged in the order 
prescribed by a map \ $\delta \ : \ \{1,  . . .  , r + t - 2\} \to 
\{0; 1, \dots, r - 2\}$. This follows inductively from the  Jacobi 
identity. Indeed, if \ $w$ \ contains a subproduct of the form \ 
$[[a, b], w_1]$ \ where \ $\pm w_1$ \ is the longest subproduct 
involving \ $\iota_{r -1}$ \ and bracketed as in 2.10, then we may 
replace \ $[[a, b], w_1]$ \ by \ $\pm [a, [b, w_1] \pm [b, [a, 
w_1]]$. 
We repeat such substitutions until \ $\iota_{r -1}$ \ occurs only 
in {\it longer} \ subproducts \ $[\iota_\ell, w_1]$ \ which are 
bracketed as in 2.10.

Now we can evaluate our Hopf homomorphisms.

\proclaim{Proposition 2.11}
Assume \ $n \ge 2$. 
Given any integers \ $s, t \ge 0$, let \ $\gamma \in \sum_{r + s - 
2}$ \ be a permutation and let \ $\delta \ : \{1, \dots, r + t - 
2\} \to \{0; 1, \dots, r -2\}$ \ be a map such that \ 
$\delta^{- 1} (\{j\})$ \ consists of precisely one element for \ 
$j = 1, \dots, r -2$.

If \ $s = t$ \ and \ $\delta$ \ is the \lq\lq contraction\rq\rq\ 
of \ $\gamma$ {\rm(}i.e. $\delta (i) = \max (\gamma (i) - t, 0)$ \ for 
$i = 1, \dots, r + t - 2)$, then for all \ $k \in \Bbb Z$ \ the 
composite map
$$
\pi_k (S^{tn + |q| - r - t + 2}) @>{\iota_{\delta *}}>> \pi_k (W) 
@>{H_{s, \gamma}}>> \pi^S_{k - sn - |q| + r + s - 2}
$$
(cf.\ {\rm 2.6}\ and {\rm 2.10}) equals the stable suspension \ $E^{\infty}$ \ 
up to a \ $\pm$--sign.

In all other cases \ $H_{s, \gamma} \circ \iota_{\delta *} \equiv 
0$.

Moreover, for \ $s \ge 0$ \ \ $H_s$ \ vanishes on all those Hilton 
summands of \ $\pi_* (W)$ \ which are not contained in \ $\pi''_* 
(W, s)$.
\endproclaim

\demo{Proof}
If \ $s \le t$, \ the first two claims follow from
arguments (based e.g.\ on the fibrewise 
intersection approach to Hopf invariants) similar to those in the 
proof of theorem 3.1 in \cite{9}. Pinching (or, 
equivalently, replacing certain submanifolds by several \lq\lq 
parallel\rq\rq\ copies) adds no difficulties since all 
codimensions are strictly larger than $1$. 

If \ $s > t$, \ then the intersection approach to \ $H_{s, \gamma} 
(\iota_\delta)$ \ requires two steps. At first we carry out the 
iterated intersection procedure described in \cite{9}, 
pp.\ 306--307, until we are confronted with the overcrossing locus
$$
g | N_{\gamma (1)} \ > \dots > \  g | N_{\gamma (s -t)} \ \ \ 
> \bar g | \Cal I
$$
as in \cite{9}, (12). Here \ $\Cal I$, \ if non empty, 
consists of a single point and hence is not nulbordant. But \ 
$N_{\gamma (1)}$ \ has a standard nulbordism \ $B$; \ thus we may 
replace overcrossings again by intersections, but in the reverse 
order, i.e.\ starting with \ $B \cap N_{\gamma (2)}$ \ etc. All 
necessary intermediate nulbordisms exist (since no \ $N_{\gamma 
(i)}$ \ corresponds to \ $\iota_{r -1}$) and they have strictly 
positive codimensions so that they miss \ $\Cal I$ \ generically. 
As in \cite{9}, top of p.\ 309, this argument can be 
applied fibre by fibre to show that \ $H_{s, \gamma} \circ 
\iota_{\delta_*} \equiv 0$. 

Finally, if a Whitehead product \ 
$\iota_\delta$ \ as in 2.10 involves a factor \ $\iota_j$ \ for 
some \ $j, \ 1 \le j \le r -2$, \ more than once then \ $H_{s, 
\gamma} \circ \iota_{\delta^*} \equiv 0$ \ since in the iterated 
(standard) intersection process \ $H_{s, \gamma}$ \ can register 
only the \lq\lq innermost\rq\rq\ factor \ $\iota_j$ \ and treats 
all outer ones as if they were zero. With a little extra care such 
an argument can also be adapted to \ $j = r -1$. In view of the 
discussion following 2.10 we conclude that \ $H_{s, \gamma} \circ 
w_* \equiv 0$ \ for every basic Whitehead product \ $w \in \pi_* 
(W)$ \ which involves some factor \ $\iota_j, \ 1 \le j \le r -1$, 
more than once (or not at all) or which involves \ $\iota_0 \ \ t$ 
times, $t \ne s$.  \hfill $\blacksquare$
\enddemo

\proclaim{Theorem 2.12}
Assume \ $n \ge 2$. Given any integers \ $k$ \ and \ $s \ge 0$, 
the {\rm(}restricted{\rm)}\ Hopf homomorphism \ $H_s| = \ H_s | 
\pi''_k (W; s)$ \ {\rm(}cf.\ {\rm2.7}\ and {\rm2.9)}\ fits into 
the commuting diagram
%
$$
\BoxedEPSF{diag212.eps} 
$$
%
where \ $\sum w_{\ell *}$ \ is the Hilton isomorphism, \ 
$E^{\infty}$ \ denotes stable suspension and \ $D_{s *}$ \ is 
defined by the operation of an invertible matrix \  $D_s \in G L 
(u_{r, s}; \Bbb Z)$.

In particular, if \ $k \le 2 (|q| - r + 1 + s (n -1))$, then \ $H_s 
|$ \ is bijective.

If even \ $k \le 2 (|q| - r + 1)$, then \ $H$ \ {\rm(}cf.\ {\rm 2.8)} 
restricts to the isomorphism {\rm(}cf.\ {\rm 2.9)}
$$
H | \ : \ \pi''_k (W) \ @>{ \ \ \ \cong \ \ \ }>> \ \bigoplus_{s 
\ge 0} \left( \pi^S_{k - |q| + r - 2 - s (n -1)}\right)^{u_{r, s}} \ .
$$
\endproclaim

\demo{Proof}
In view of \cite{6}, p.\ 155, the number of basic Whitehead 
products \ $w_{\ell}$ \ generating \ $\pi''_k (W; s)$ \ turns out to 
coincide with \ $u_{r, s}$ (compare 2.5). Up to transposition and 
\ $\pm$ signs, the \ $u_{r, s} \times 
u_{r, s}$--matrix \ $D_s$ \ consists of the integer coefficients 
encountered when we express these \ $w_{\ell}$ \ in terms of the 
elements \ $\iota_\delta$ \ as in 2.10 (for \ $t = s$) or, 
equivalently, of their \ $\pi''_* (W, s)$--components.
But these, in turn, are \ $\Bbb Z$--linear 
combinations of the basic Whitehead products \ $w_{\ell}$. Thus \ $D_s$ 
\ is invertible.

As in 2.11 we evaluate the Hopf invariants \ $H_s$ \ fibrewise via the 
intersection approach (compare \cite{9}, p.\ 308--309); this 
yields our commutativity claim. \hfill $\blacksquare$
\enddemo

\vskip2mm
It remains to discuss {\it the case when} \ $n = 1$.

Let \ $\widetilde W$ \ denote the universal covering space of \ 
$W$. We may think of it as a real line \ $\Bbb R$, with a wedge \ 
$\widetilde W_g = \bigvee^{r -1}_{j = 1} S^{q_j}_{j, g}$ \ 
attached at every integer \ $g \in \Bbb Z \subset \Bbb R$. We 
write
$$
c \ : \ \widetilde W \ \longrightarrow \ \Bbb R
\tag 2.13
$$
for the obvious \lq\lq level\rq\rq\ map. It is compatible with the 
\lq\lq shifts\rq\rq\ (deck transformations) by elements in \ $\Bbb 
Z = \pi_1 (W) = \pi_1 (S^1)$.

On the other hand, there is a canonical isomorphism \ $\pi_* 
(\widetilde W) \cong \pi_* (\vee_g \widetilde W_g)$ \ and we may 
choose a Hilton decomposition based on Whitehead products which 
involve the inclusions \ $\iota_{j, g} : S^{q_j}_{j, g} \subset 
\vee_g \widetilde W_g$ \ (where $1 \le j \le r - 1$ and $g \in \Bbb Z)$.

This applies also to our original wedge \ $W$ \ via the 
isomorphism
$$
p_* \ : \ \pi'_* (\widetilde W) \ @>{ \ \ \cong \ \ }>> \ \pi'_*  
(W)
\tag 2.14
$$
of {\it reduced}\ homotopy groups (compare 2.3 and also \cite{10}, 2.6 and 2.10) where \ $p$ \ denotes the covering projection. Then \ $\pi'_* (W)$ \ is the 
direct sum of all Hilton summands corresponding to basic Whitehead 
products which for each \ $j = 1, \dots, r -1$ \ involve at least 
one factor of the form \ $\iota_{j, g}, g \in \Bbb Z$.
Define the subgroup
$$
\pi''_* (W) \ \subset \ \pi'_* (W)
\tag 2.15
$$
to be the direct sum of all such summands which involve {\it 
precisely} \ one such factor for each $j$.

Now let us evaluate our Hopf homomorphisms on an arbitrary class \ 
$[v] \ \in \ \pi'_k (W)$. We adopt the intersection approach 
outlined in \cite{9}, pp.\ 306--307. 
\vskip0.2cm
%
$$
\EPSFxsize=12.3cm 
\BoxedEPSF{umgebung.eps} 
$$
%
\vskip-0.5cm
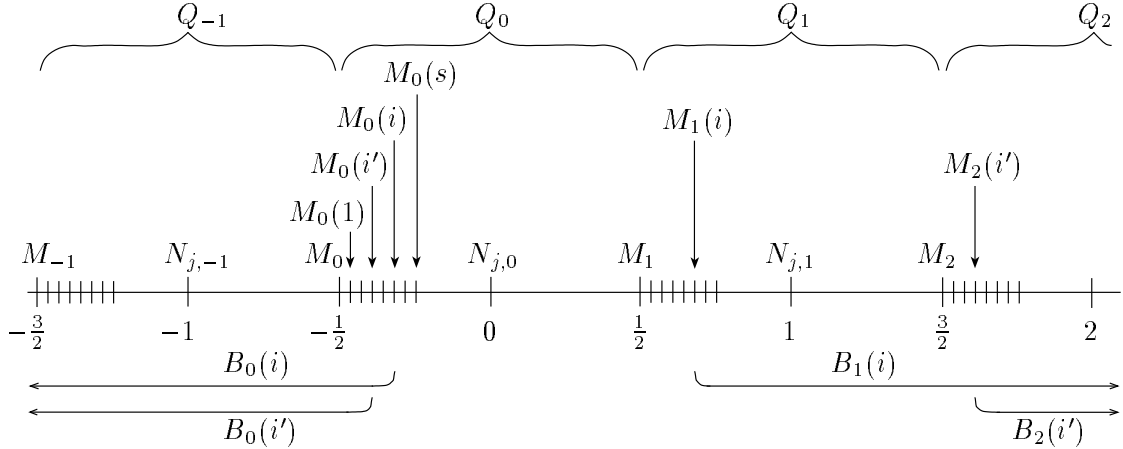
\botcaption{Figure 2.16} The arrangement of the manifolds \ $N_{j, 
g}$ and $M_g(i)$ \newline
(Here link components and other subsets of $\Bbb R^k$ are listed, 
together \newline  with their images in \ $\Bbb R$ \ under \ $c \circ 
\widetilde v$)
\endcaption

\vskip0.7cm
Let \ $\widetilde v 
\ : \ S^k \to \widetilde W$ \ be the lifting of a generic 
representative $v$. Then for \ $j = 1, \dots, r -1$ \ the link 
component
$$
N_j \ = \ v^{- 1} (\{z_j\}) \ \subset \ \Bbb R^k \ = \ S^k - \{*\} 
$$
corresponding to \ $S^{q_j}_j$ \ is the disjoint union of the 
manifolds
$$
N_{j, g} \ = \ \widetilde v^{-1} (\{g \widetilde z_j\}) , \ \ g \in 
\Bbb Z \ ;
\tag 2.17
$$
here \ $z_j \ \in \ S^{q_j}_j - \{*\}$ \ and \ $g \widetilde z_j \ 
\in \ S^{q_j}_{j, g} - \{*\}$ \ are regular values of \ $v$ and 
$\widetilde v$. 
Similarly the link component \ $M \ = \ v^{- 1} (\{ z_0\}), \ \ z_0 
= (- 1, 0) \in S^1 - \{*\} \ \subset \ W$, is the disjoint union of 
the hypermanifolds \ $M_g \ = \ \widetilde v^{- 1} (\{g - 
\frac12\}), \ g \in \Bbb Z$, which decompose \ $\Bbb R^k$ \ into 
the strips \ $Q_g = (c \widetilde v)^{-1} [g - \frac12, g + 
\frac12]$. Note that \ $Q_0$ is unbounded and that
$$
N_{j, g} \ = \ Q_g \cap N_j \ \ \ \text{for}\ j = 1, \dots, r - 1 \ 
\text{and all} \ g \in \Bbb Z \ .
$$
For \ $s \ge 0$ \ pinching \ $S^1$ \ $s$ times amounts to 
replacing \ $M$ (and $M_g$) by $s$ nearby \lq\lq parallel\rq\rq\ 
copies \ $M (i)$ \ (and $M_g (i)), \ i = 1, \dots, s$, (where e.g.
$$
M_g (i) \ = \ (c \widetilde v)^{- 1} (\{ g - \frac12 + 
i \varepsilon \}) \ , \ \ g \in \Bbb Z \ ,
\tag 2.18
$$
for some small fixed \ $\varepsilon > 0$, compare figure 2.16).

Given an $s$--monotone permutation \ $\gamma \in \sum_{r, s}$ \ 
(cf.\ 2.4), it will be convenient to calculate the Hopf invariant 
\ $H_{s, \gamma} [v]$ \ via iterated intersections as in \cite{9}, (12), but in the {\it reverse order}\ (this changes the 
value only by a fixed sign which is obvious from the overcrossing 
description of \ $H_{s, \gamma} [v]$, \ cf. p. 306 of \cite{9}). Thus if \ $N (\gamma (1)), N (\gamma (2)), \dots, N 
(\gamma (r + s - 2)), \ N_{r - 1}$ \ are the link com\-po\-nents of \linebreak
($\pinch_s \vee \id)_* [v]$, listed in the order given by \ 
$\gamma$, intersect a framed (singular) nulbordism of \ $N (\gamma 
(1))$ in $\Bbb R^k$ \ with \ $N (\gamma (2))$ \ to obtain a closed 
framed singular manifold \ $N (\gamma (1), \gamma (2))$; similarly, 
intersect a nulbordism of \ $N (\gamma (1), \gamma (2))$ \ in \ 
$\Bbb R^k$ \ with \ $N (\gamma (3))$, etc., until finally a 
nulbordism of \ $N (\gamma (1), \dots, \gamma (r + s - 2))$ is 
intersected with \ $N_{r - 1}$ \ to yield \ $\pm H_{s, \gamma} 
[v]$.

Since by assumption \ $[v]$ \ lies in the reduced homotopy group \ 
$\pi'_k (W)$, all intermediate intersections which occur in this 
process are  indeed nulbordant. In fact, often we have even 
canonical nulbordisms. E.g. $M_g (i)$ \ bounds the compact 
manifold \ $B_g (i)$ \ defined to be \ $(c \widetilde v)^{- 1} [g 
- \frac12 + i \varepsilon, + \infty)$ \ if \ $g \ge 1$, or, with 
opposite sign, \ $(c \widetilde v)^{- 1} (- \infty, g - \frac12 + 
i \varepsilon]$ \ if $g \le 0$. Thus \ $B (i) = \bigcup_{g \in 
\Bbb Z} \ B_g (i)$ \ is a nulbordism of \ $M (i)$ \ which, for 
every \ $g \in \Bbb Z$, covers most of the strip \ $Q_g$ -- and 
hence all of \ $N_{j, g}, \ j = 1, \dots, r -1$, -- precisely $g$
times. Figure 2.16 may help to visualize such statements.

Similarly, \ $B_g (i) \cap B_{g'} (i')$ \ is a nulbordism of \ 
$\pm B_g (i) \cap M_{g'} (i')$ \ if \ $1 < g + i \varepsilon < g' 
+ i'\varepsilon$ \ or \ $1 > g + i \varepsilon > g' + i' 
\varepsilon$ \ etc.

\example{Example 2.19: {\sl\bf r} = 2}
Here the generic base point preserving map
$$
\widetilde v \ : \ S^k \ \longrightarrow \ \widetilde W \sim 
\bigvee_{g \in \Bbb Z} S^{q_1}_{1, g}
$$
corresponds to the (finite) framed link
$$
N_1 \ = \coprod_{g \in \Bbb Z} N_{1, g} \ \subset \ \Bbb R^k = S^k 
- \{*\}
$$
For each \ $g \in \Bbb Z$ \ the (stable) framed bordism class
$$
[N_{1, g}] \ \in \ \Omega^{fr}_{k - q_1} \ \cong \ \pi^S_{k -q_1}
$$
is given by the stable suspension of the composite of \ 
$\widetilde v$ \ with the map which collapses all wedge summands \ 
$S^{q_1}_{1, g'}, \ g' \ne g$, \ to a point.

For each \ $s \ge 0$ \ there is only one $s$--monotone permutation 
\ $\gamma \in \sum_{2, s}$. We can apply the previous discussion 
to the resulting Hopf homomorphism \ $H_s = H_{s, \gamma}$ \ and 
obtain
$$
\pm H_s ([v]) \ = \sum_{g \in \Bbb Z} {g + s - 1 \choose s} [N_{1, 
g}] \ \in \pi^S_{k -q_1} \ \ \ .
\tag 2.20
$$
Indeed, \ $H_0 ([v])$ \ does not involve any of the hypermanifolds 
\ $N_g (i)$ \ and \ $\pm H_0 ([v])$ \ is just equal to \ $[N_1]´= 
\sum [N_{1, g}]$. Next,
$$
\pm H_1 ([v]) \ = \ [B (1) \cap N_1] \ = \ \sum_{g'} [B_{g'} (1) \cap 
N_1]
$$
counts each link component \ $N_{1, g}$ \ precisely \ $g$ \ times 
since it is contained in just so many $0$--codimensional bordisms 
\ $B_{g'} (1)$ \ (which have the opposite orientations for \ $g' 
\le 0$ \ so that the signs are also correct). In the same spirit 
but more generally, let \ $b_{s, g}$ \ denote the integer which 
counts the relevant $s$--fold intersections \ $B_{g_1}(1) \cap \dots 
\cap B_{g_s} (s)$ \ containing \ $N_{1, g}$, \ $s \ge 1, \ g \in 
\Bbb Z$. (By \lq\lq relevant\rq\rq\ we mean that the intersection 
can contribute to \ $H_s [v]$, i.e. $B_{g_i} (i) \supset M_{g_{i + 
1}} (i + 1)$ \ for $1 \le i < s)$. Comparing such intersections to those 
which contain \ $N_{1, g - 1}$, \ we see that
$$
b_{s, g} \ = \ b_{s, g -1} \ + \ b_{s -1, g}
\tag 2.21
$$
(figure 2.16 may again be helpful here). The identity
$$
b_{s, g} \ = \ {g + s - 1 \choose s}
\tag 2.22
$$
and formula 2.20 follow now by induction over $s$ (and, for fixed 
$s$, over $|g|$). \hfill $\blacksquare$
\endexample

The next result shows that the sequence \ $\{[N_{1, g}]\}_{g \in 
\Bbb Z}$ \ is completely determined by the sequence \ $\{H_s 
([v])\}_{s \ge 0}$ \ of Hopf invariants.

\proclaim{Lemma 2.23}
For any abelian group A the homomorphism
$$
d \ = \ \prod_{s \ge 0} d_s \ : \ \bigoplus_{g \in \Bbb Z} A \ 
\longrightarrow \ \prod_{s \ge 0} A \ \ ,
$$
defined by \ $d_s ((a_g)_{g \in \Bbb Z}) \ = \ \sum_{g \in \Bbb Z} 
{g + s - 1 \choose s} a_g$, \ is injective.
\endproclaim

\demo{Proof}
Given any integers \ $n \ge 0$ \ and \ $n_0$, consider the \ $(n + 
1) \times (n + 1)$--matrix \ $M (n, n_0)$ \ with entries
$$
b_{s, g} \ = \ {g + s - 1 \choose s} , \  \ 0 \le s \le n, \ \ \ 
n_0 - n \le g \le n_0 \ .
$$
When we subtract the \ $n^{\text{th}}$ column from the \ $(n+ 
1)^{\text{st}}$ one, the \ $(n - 1)^{\text{st}}$ column from the \ 
$n^{\text{th}}$ column and so forth, the top row takes the form \ 
$(1, 0, \dots, 0)$ \ and -- in view of 2.21 -- the \ $n \times 
n$--submatrix in the lower right hand corner coincides with \ $M 
(n - 1, n_0)$. Thus
$$
\det M (n, n_0) \ = \ \det M (n - 1, n_0) \ = \dots = \ \det M (0, 
n_0) \ = 1 \ .
$$
Since every element of \ $\ker d$ \ lies already in the kernel of 
the endomorphism \ $M (n, n_0)$ \ on \ $\bigoplus_{n_0 - n \le g 
\le n_0} A$ \ for suitable \ $n$ and $n_0$, our lemma 
follows. \hfill $\blacksquare$
\enddemo

Next let \ $r \ge 2$ \ be arbitrary. Given an \ $(r - 1)$--tuple \ 
$(g) = (g_1, \dots, g_{r - 1}) \in \Bbb Z^{r - 1}$ of integers and 
a permutation \ $\bar\gamma \in \sum_{r - 2}$, let
$$
h_{(g), \bar\gamma} \ : \ \pi'_k (W) \ \longrightarrow \ \pi^S_{k 
- |q| + r - 2}
\tag 2.24
$$
be the homomorphism which maps a homotopy class \ $[v] \in \pi'_k 
(W)$ or, equivalently, the corresponding framed link \ $\amalg \
N_{j, g} \subset \Bbb R^k$ \ as in 2.17, to the Hopf invariant 
(cf.\ \cite{9}, \S\ 3)
$$
h_{\bar\gamma} \ \left(N_{1, g_1} \ \amalg \ \dots \ \amalg  \ 
N_{r - 1, g_{r - 1}} \right)
$$
of the indicated sublink (on the homotopy level the selection of 
this sublink is induced by the obvious map which collapses \ 
$\widetilde W$ \ to \ $\bigvee^{r - 1}_{j = 1} S^{q_j}_{j, g_j}$, 
compare the discussion preceding 2.14).

Now, for \ $s \ge 0$, \  consider an \ $s$--monotone permutation 
\ $\gamma \in \sum_{r, s}$. Then the list of values of \ $\gamma$  
\ takes the form
$$
(1, \dots, s_1; \bar\gamma (1) + s; s_1 + 1, \dots, s_1 + s_2; 
\bar\gamma (2) + s; \dots; \bar\gamma (r - 2) + s; s_1 + \dots + 
s_{r -2} + 1, \dots, s)
\tag 2.25
$$
where \ $(s) = (s_1, \dots, s_{r -1}) \in \Bbb Z_+^{r - 1}$ \ is 
an \ $(r - 1)$--tuple of nonnegative integers whose sum equals \ $ 
s$, and $\bar\gamma \in \sum_{r -2}$; in words: \ $s_i$ (strictly 
increasing) values lying in \ $\{1, \dots, s\}, \ i = 1, \dots, r 
-1$, alternate with the values \ $\bar\gamma (i) + s \ \in \{s + 
1, \dots, r + s - 2\}$ \ which appear in the order given by the 
permutation \ $\bar\gamma$. We obtain a bijection between such 
pairs \ $((s), \bar\gamma)$ \ and the corresponding \ 
$s$--monotone permutations \ $\gamma ((s), \bar\gamma) \ \in \ 
\sum_{r, s}$ \ defined by 2.25.

\proclaim{Proposition 2.26}
For \ $\gamma = \gamma ((s), \bar\gamma)$ \ as above and for every 
homotopy class \ $[v] \ \in \ \pi'_k (W),  \ k \in \Bbb Z,$ \ 
{\rm(}cf.\ {\rm2.3)} we have
$$
H_{s, \gamma} ([v]) \ = \ \varepsilon \cdot \sum_{(g) = (g_{1}, 
\dots, g_{r -1}) \in \Bbb Z^{r - 1}} \ \prod^{r -1}_{j = 1} {\bar 
g_j - \bar g_{j - 1} + s_j - 1 \choose s_j} \ h_{(g), \bar\gamma} 
\ ([v])
$$
where \ $\varepsilon = \pm 1$ \ is a fixed sign depending only on 
\ $\gamma$ \ and \ $k$, \
and \ $\bar g_j := g_{\bar\gamma (j)}$ \ for $j = 1, \dots, r -2$ 
(and \ $\bar g_0 := 0$ \ and \ $\bar g_{r - 1} := g_{r - 1})$.
\endproclaim

\demo{Proof}
Adopting the intersection approach (see the discussion which 
follows 2.18) we iterate the procedure used in example 2.19 
(compare 2.20). 
After the first step (i.e.\ after the \ $s_1$--fold intersection 
with \lq\lq standard\rq\rq\ $0$--codimensional nulbordisms), each 
link component \ $N_{\bar\gamma (1), \bar g_1}$ (cf.\ 2.18), $\bar 
g_1 \in \Bbb Z$, is counted with the multiplicity \ 
$b_{s_1, \bar g_1}$ (cf.\ 2.22).

If \ $r > 2$ \ then each \ $N_{\bar\gamma (1), \bar g_1}$, in turn, 
allows a (generic, singular) nulbordism $B$ in $\Bbb R^k$ 
(use 2.15 and the fact that \ $[v]$ \ is reduced). 
Since \ $c \circ \widetilde v$ \ maps \ 
$\partial B$ \ to \ $\bar g_1$ \ (cf.\ figure 2.16), $B \cap M_g 
(i)$ \ (cf.\ 2.18), \ $g \in \Bbb Z$, bounds the intersection of \ 
$B$ \ with \ $(c \widetilde v)^{-1} [g - \frac12 + i \varepsilon, 
\infty)$ \ if \ $g > \bar g_1$, and with \ $(c \widetilde v)^{-1} 
(- \infty, g - \frac12 + i \varepsilon]$ \ if $g \le \bar g_1$. In 
the second step of our iteration we have to intersect \ $B$ \ 
$s_2$--times in this fashion and then with \ $N_{\bar\gamma (2)}$. 
This results in \ $B \cap N_{\bar\gamma (2), \bar g_2}, \ \bar g_2 
\in \Bbb Z$, being counted with the (additional) multiplicity \ 
$b_{s_2, \bar g_2 - \bar g_1}$ \ (use the same arguments as in 
example 2.19, but with \ $g$ \ replaced by \ $g - \bar g_1$).

Continuing this process and intersecting with \ $N_{r - 1}$ \ in 
the last step, we obtain the indicated linear combination of the 
Hopf invariants \ $h_{(g), \bar\gamma} [v]$. Each of these can be 
evaluated via the intersection approach since \ $[v]$ \ is reduced
and hence all (intermediate intersections involving) sublinks with
strictly less than \ $r - 1$ \ components \ $N_{j, g}$ \ allow
nulbordisms (compare 2.14). \hfill $\blacksquare$
\enddemo

Now we are ready to establish the analogue of theorem 2.12.

\proclaim{Theorem 2.27}
Assume \ $n = 1$. Then for all integers \ $k$ \ the Hopf 
homomorphism \ $H$ {\rm(}cf.\ {\rm 2.8)}, when restricted to \ $\pi''_k (W) $ 
\ {\rm(}cf.\ {\rm2.15)}, fits into the commuting diagram
%
$$ 
\BoxedEPSF{diag227.eps} 
$$ 
%
where \ $\sum w_{\ell*}$ \ is the Hilton isomorphism {\rm(}cf.\ {\rm2.15)}, 
\ $E^{\infty}$ \ denotes stable suspension and \ $D$ \ is a group 
monomorphism. Moreover, \ $H$ \ vanishes on all those Hilton 
summands of \ $\pi_k (W)$ \ which are not contained in \ $\pi''_k 
(W)$.

In particular, if \ $k \le 2 (|q| - r + 1)$, then \ $H$ \ is 
injective on \ $\pi''_k (W)$.
\endproclaim

\demo{Proof}
According to the last proposition there is a commuting diagram
%
$$ 
\BoxedEPSF{diag227a.eps} 
$$ 
%
\flushpar
where the range of \ $H$ \ is indexed by all \ $s$--monotone 
permutations, \ $s \ge 0$, in the form \ $\gamma = \gamma ((s), 
\bar\gamma)$ \ (compare 2.25) and \ $D'$ \ is defined by the \ 
$\Bbb Z$--linear combinations in 2.26.

First we want to prove that \ $D'$ -- or, equivalently, \ 
$D'_{\bar\gamma}$ \ for each \ $\bar\gamma \in \sum_{r -2}$ -- is 
a monomorphism. After reparametrizing suitably we must show this only 
for the map \ $D'_{r -1}$ \ defined by
$$
D'_{r -1} ((a_{(g)})_{(g) \in \Bbb Z^{r -1}}) = \left( \sum_{(g)} 
\prod^{r -1}_{j = 1} {g_j + s_j - 1 \choose s_j} 
a_{(g)}\right)_{(s) \in \Bbb Z^{r -1}_+} \ .
$$
Thus suppose that
$$
a \ = \ (a_{(g)})_{(g) = (g_1, \dots, g_{r -1}) \in \Bbb Z^{r -1}}
$$
lies in the kernel of \ $D'_{r -1}$. Then for every \ $(s) = (s_1, 
\dots, s_{r -1}) \in \Bbb Z^{r -1}_+$ \ the expression
$$
\sum_{g_{r -1} \in \Bbb Z} {g_{r -1} + s_{r -1} - 1 \choose s_{r 
-1}} \left( \sum_{(g_1, \dots, g_{r -2}) \in \Bbb Z^{r -2}}  \ 
\prod^{r -2}_{j =1} {g_j + s_j -1 \choose s_j} a_{(g_1, \dots, g_{r 
- 2}, g_{r -1})} \right)
$$
vanishes, and so does the sum to the right hand side, for every 
fixed \ $g_{r -1} \in \Bbb Z$, by lemma 2.23. Note that such sums, 
indexed by \ $(s_1, \dots, s_{r -2})$, constitute the value of \ 
$a_{(-, \dots, -, g_{r -1})}$ under \ $D'_{r -2}$. Thus we can 
iterate our argument until the injectivity of \ $D'_1 = d$ (cf.\ 
2.23) implies that \ $a = 0$.

Next we have to compare, for fixed \ $(g) \in \Bbb Z^{r -1}$, the 
values of \ $h_{(g), \bar\gamma}, \ \bar\gamma \in \sum_{r -2}$, to 
the stable suspensions of the Hilton components which correspond 
to the basic Whitehead products in \ $\iota_{1, g_1}, \dots, 
\iota_{r -1, g_{r -1}}$. When we deal with \ $\pi''_k (W)$, each 
of these factors is involved precisely once, and the comparison 
can be expressed via an invertible matrix \ $D_0$ as in theorem 
2.12. Composing \ $\bigoplus_{(g)} D_{0 *}$ \ with \ $D'$ \ we 
obtain the desired monomorphism \ $D$ \ in 2.27. 

Finally, consider a basic Whitehead product which, for some \ $1 
\le j \le r -1$, involves more than one factor of the form \ 
$\iota_{j, g}, \ g \in \Bbb Z$ \ (or none at all). As in the proof 
of proposition 2.11 we conclude that \ $H \circ w_* \equiv 0$. 
Clearly, also \ $H | \pi_1 (W) \equiv 0$. \hfill $\blacksquare$
\enddemo

In view of theorems 2.12 and 2.27 it is of interest to compare the 
subgroups \ $\pi'_k (W)$ \ and \ $\pi''_k (W)$ of $\pi_k (W), \ k 
\in \Bbb Z$ \ (cf.\ 2.3, 2.9, and 2.15).

\proclaim{Lemma 2.28}
Assume \ $n \ge 1$. If
$$
k \ \le \ |q| \ + q_j - r \ \ \ \text{for} \ j = 1, \dots, r -1 \ 
,
$$
then \ $\pi'_k (W) = \pi''_k (W)$.
\endproclaim

Indeed, \ $\pi''_k (W)$ \ {\rm(}or $\pi'_k (W)$, resp.{\rm)} \ is 
the direct sum of the Hilton summands corresponding to those basic 
Whitehead products which for every \ $j = 1, \dots, r -1$ \ 
involve the inclusion \ $\iota_j \ : \ S^{q_j} \subset W$ \ 
{\rm(}if \ $n \ge 2${\rm)} \ or an inclusion of the form \ 
$\iota_{j, g} \ : \ S^{q_j} \subset \widetilde W$ \ {\rm(}if \ $n 
= 1)$ \ precisely once {\rm(}or at least once, resp.{\rm)}. In the 
indicated dimension range the additional Hilton summands in \ 
$\pi'_k (W)$ \ are all trivial.

\vskip7mm

\specialhead  3. \ \ Higher $\mu$-invariants and their 
geometry
\endspecialhead

From now on let \ $M$ \ denote the manifold \ $S^n \times \Bbb 
R^{m -n}$ \ (where \ $1 \le n < m \ge 3$) unless mentioned otherwise.
Pick a suitable base 
point \ $y^0 = (y^0_1, \dots, y^0_r)$ \ of the configuration space
$$
\widetilde C_r (M) \ = \left\{ (y_1, \dots, y_r) \in M^r \ | \ y_i \ne 
y_j \ \ \text{for} \ 1 \le i \ne j \le r\right\}
\tag 3.1
$$

Suitable embeddings \; $S^n = S^n \times \{ x_0\}  \subset M$ \; and
\; $(\bigvee^{r -1} B^m_j, *)  \subset (M, y^0_r)$ \; 
(such that the base point \ $y^0_j$ \ 
lies in the interior of the ball \; $B^m_j, \ j = 1, \dots
r -1)$ \; yield a homotopy equivalence
$$
S^n \vee \bigvee^{r -1}_{j =1} S^{m -1}_j \ \hookrightarrow \
M - \{y^0_1, \dots, y^0_{r -1} \} 
$$
which is canonical up to homotopy (if in the case \ $n = m -1$ \ we also
require \ $x_0 \in \Bbb R$ \ to be so small that \ $y^0_1, \dots,
y^0_{r -1} \in S^{m -1} \times (x_0, \infty)$). Together 
with the \lq\lq fibre inclusion\rq\rq
$$
M - \{y^0_1, \dots, y^0_{r -1}\} \ \longrightarrow \ \widetilde C_r (M), 
\quad y \longrightarrow (y^0_1, \dots, y^0_{r -1}, y) \ ,
$$
and an obvious quotient map this induces the injective composite 
map
$$
%
\EPSFxsize=12.3cm 
\BoxedEPSF{diag32.eps} 
\tag 3.2
$$
(compare \cite{10}, 2.5) where
$$
|p| \ := \ p_1 + \dots + p_r \ \ \ .
\tag 3.3
$$

Now we are ready to define higher \ $\mu$-invariants. First we 
work in a setting where the base points \ $* \in S^{p_j}$ \ are 
preserved by link maps and link homotopies. Later we will also 
comment on the question as to when the resulting invariants 
remain unchanged under base point free link homotopies.

Thus let
$$
f \ = \ f_1 \amalg \dots \amalg f_r \ : \ S^{p_1} \amalg \dots 
\amalg S^{p_r} \ \longrightarrow \ S^n \times \Bbb R^{m -n}
$$
be a \ $\km$-Brunnian link map which, in addition, preserves 
base points, i.e.\ $f_j (*) = y^0_j$ \ for $j = 1, \dots, r$. Then 
according to \cite{10}, 2.4, there is a unique reduced 
homotopy class \ $\km' (f)$ as in diagram 3.2 above such that 
\ $\quot^* \circ \incl_* (\km' (f))$ \  equals the base point 
preserving version \ $\km^b (f)$ \ of the \ 
$\kappa$-invariant of \ $f$ (compare 1.2). For \ $s \ge 0$ \ 
define
$$
\mum^{(s)} (f) \ := \ H_s (\km' (f)) \in \bigoplus^{u_{r, s}} 
\pi^S_{|p| - s (n -1) - (r -1) (m -2) -1} \ 
\tag 3.4
$$
(cf.\ 2.7).

Note that \ $\mum^{(0)} (f)$ \ coincides with the invariant \ 
$\mum (f)$ \ discussed in \cite{10}, \S\ 2. Indeed, the 
homomorphism \ $c_a$ \ used there is induced by the map which 
collapses \ $S^n$; hence \ $H_0 = h \circ c_a$ \ (if \ $n = m -1$ 
\ we assume here that the arcs involved in \ $c_a$ converge to \ 
$S^{m -1} \times \{ + \infty\})$.

The main results 2.12 and 2.27 of the previous section, together 
with lemma 2.28, now imply

\proclaim{Theorem 3.5}
If \ $p_1 + \dots + p_r \le r (m -2)$, then the base point 
preserving \ $\kappa$-invariant \ $\km^b (f)$ \ of a 
$\km$-Brunnian link map \ $f$ \ into \ $M = S^n \times \Bbb R^{m - 
n}, \ 1 \le n < m \ge 3$, \ contains precisely as much information 
as the sequence \ $\{ \mum^{(s)} (f)\}_{s \ge 0}$ \ which starts 
with \ $\mum (f)$.
\endproclaim

What is the  geometric meaning of the remaining \lq\lq higher \ 
$\mu$-invariants\rq\rq\ in this sequence?

In order to answer this question, fix any orientation preserving
smooth embedding
$$
\eta \ : \ S^n \times B^{m -n} \ \hookrightarrow \ \Bbb R^m
\tag 3.6
$$
where \ $B^{m -n}$ \ is the compact \ $(m -n)$--dimensional unit ball
with interior \; $\overset\circ\to {B}^{m -n} \cong \Bbb R^{m -n}$.
Thus \ $\eta$ \ defines an inclusion of \; $S^n \times \Bbb R^{m -n}$ \
into \ $\Bbb R^m$, as well as a meridian
$$
\eta_z \ : \ S^{m - n - 1} \ = \ \{ z\} \times \partial B^{m -n} \ \
{\overset{\eta|}\to \hookrightarrow} \ \ \Bbb R^m
$$
for every \ $z \in S^n$.

Now, given any link map in \ $S^n \times \Bbb R^{m -n}$, \; add \ $s$ \
meridians \; $\eta_{z_1}, \dots, \eta_{z_s}$ \ (at pairwise distinct points
\; $z_1, \dots, z_s \in S^n$). Up to link 
homotopy the resulting augmented link map
$$
f^{(s)} \ = \ \left( \coprod^s_{j = 1} \eta_{z_j}\right)\amalg 
\eta \circ f \ \ : \ \ \coprod^s_{j =1} S^{m -n -1} \amalg \coprod^r_{i = 1}
S^{p_i} \ @>>{ \qquad }> \ \Bbb R^m \qquad ,
\tag 3.7
$$
$s \ge 0$, depends only on the link homotopy class of \ $f$ in $M$.

\proclaim{Theorem 3.8}
Assume \; $|p| \le r (m -2)$ \ and \ $s \ge 0$. Let \; $f : \amalg S^{p_i}
\longrightarrow M = S^n \times \Bbb R^{m -n}$ \ be a $\km$--Brunnian and
base point preserving link map.

If the augmented link map \; $f^{(s)}$ in $\Bbb R^m$ \; 
{\rm(}cf.\ {\rm 3.7)}\ is
\ $\kappa_{\sssize {\Bbb R^m}}$--Brunnian {\rm(}or, equivalently,
\; $\mum^{(t)} (f) = 0$ \ for all $t < s$ {\rm)}, 
then for all \; $\gamma \in
\Sigma_{s + r -2}$
$$
\mu_{\sssize M, \gamma}^{(s)} (f) \ := \ H_{s, \gamma} (\km' (f)) =
\mu_{\sssize{\Bbb R^m},\gamma} (f^{(s)})
$$
(at least up to a fixed sign), and hence, in particular,
\; $\mum^{(s)} (f)$ \; depends only on the base point free link
homotopy class of \ $f$ in $M$.
\endproclaim

\demo{Proof}
Since \ $f$ \ is \ $\km$--Brunnian, the product map
\ $\widehat f$ \ (cf.\ 1.2)  can be deformed in a base point 
preserving fashion until it factors through a map of the
form
$$
F^{(s)}_{s + r} \ : \ S^{|p|} \ \longrightarrow \ S^n \vee 
\bigvee^{r -1} S^{m -1} \ \subset \ M - \{ y^0_1, \dots,
y^0_{r -1} \} \ ({\underset{\eta}\to{\hookrightarrow}} \Bbb R^m).
$$
This extends to the new link map
$$
F^{(s)} \ = \ \amalg \ F^{(s)}_i \ : \ \coprod^s_{j =1} \ S^{m-n-1} \ 
\amalg \coprod^{r -1}_{i =1} S^0 \ \amalg \ S^{|p|} \ \longrightarrow \
\Bbb R^m
\tag 3.9
$$
which involves also the meridians \; $F^{(s)}_j =  \eta_{z_j}, \ 
j =1, \dots, s$ \ (compare 3.7), and \ $F^{(s)}_{s +i}$ \ 
which map the base point $*$ of \ $S^0$ \ to a suitable value in \ 
$\eta (S^n \times \partial D^{m -n})$ \ and the other point of \ 
$S^0$ \ to \ $y^0_i, \ i = 1, \dots, r -1$. (Recall that 
0-dimensional link components present no problems in \ $\Bbb R^m$, 
cf.\ \cite{9}). If \ $F^{(s)}$ \ is 
$\krm$-Brunnian then so is \ $f^{(s)}$, and the \ 
$\krm'$-invariants of \ $F^{(s)}$ and $f^{(s)}$ \ coincide (they 
vanish if and only if \ $F^{(s +1)}$ or, equivalently, \ 
$f^{(s +1)}$ \ is \ $\krm$-Brunnian); on the other hand, since \ 
$F^{(s)}$ is link homotopic to \ $e_* ((\pinch_s \vee \id) 
\circ F^{(s)}_{s +r})$, it follows from theorem 6.1 in 
\cite{9} that for all \ $\gamma \in \sum_{s + r - 2}$
$$
h_\gamma (\krm' (F^{(s)})) \ = \ \pm  h_\gamma ((\pinch_s \vee \id) 
\circ F^{(s)}_{s +r}) \ .
$$
But these two expressions are \ $\mu_{{\sssize\Bbb R^m}, \gamma} 
(f^{(s)}$ \ and \
$H_{s, \gamma} (\km' (f))$, resp. We conclude, in particular, 
that \ $\mu^{(s)} (f)$ \ vanishes if and only if \ $\krm' (f^{(s)})$ 
\ does (use (24) in \cite{9}). An induction over \ $s$ \ 
now completes the proof. \hfill $\blacksquare$
\enddemo

Recall that the base point free homotopy class \ $\km (f)$ (cf.\ 
1.2) is trivial precisely if its base point preserving analogue \ 
$\km^b (f)$ \ is.

\proclaim{Corollary 3.10}
Assume \ $p_1, \dots, p_r \le m -2$.

Then for all link maps \; $f : \amalg \ S^{p_i} \longrightarrow M =
S^n \times \Bbb R^{m -n}$ \; the following two conditions are
equivalent:
\roster
\item"(i)" $\km (f)$ \ is trivial {\rm(}cf.\ {\rm 1.2)} \qquad \qquad ; and
\item"(ii)" the component maps \ $f_i : S^{p_i} \longrightarrow M \sim S^n$
are nulhomotopic, \ $i = 1, \dots, r$, {\rm and} \ the invariant \ $\krm 
(f^{(s)})$ \ of the augmented link
map \ $f^{(s)}$ is trivial for every \ $s \ge 0$.
\endroster
\endproclaim

\demo{Proof}
Recall from \cite{9}, theorem 4.2, that the second part of condition
(ii) holds if and only if all (the consecutively defined) \ 
$\mur$--invariants of all sub-link maps of \ $f^{(s)}, \ s \ge 0$, vanish.

Clearly, (i) implies (ii). Assume that we have proved the converse
inductively for link maps with up to \ $r -1$ \ components. Then (ii)
implies that \ $f$ \ is \ $\km$--Brunnian and, by theorem 3.8,
even that \ $\km (f)$ is trivial (use induction over \ $s$ \ and 
theorem 3.5) \hfill $\blacksquare$
\enddemo

\remark{Remark {\rm 3.11}}
If \ $n \ge 2$, then \ $M$ is $1$--connected; thus \ $\km^b (f)$ \ 
(and hence all \ $\mum^{(s)} (f)$, when defined) are invariant even 
under {\it base point free}\ link homotopies of \ $f$ \ (compare 
\cite{10}, \S\ 3).
Moreover, the higher \ $\mu$-invariants \ 
\ $\mum^{(s)} (f), \ s \ge 1$, are in general truly new and not
already determined by the invariant \ $\mum^{(0)} (f) = \mum (f) =
\widetilde\mu_{\sssize M} (f)$ \ studied in \cite{10}. This is 
illustrated e.g.\ by example 5.9 of \cite{10} where \ $f^{(1)}$ 
\ is a higher dimensional version of the classical Borromean link 
and \ $\mum (f) = 0$ \ but \ $\mu^{(1)} (f) \ne 0$.

If \ $n = 1$ \ and \ $p_1, \dots, p_r \le m -2$, the base point 
preserving \ $\kappa$--invariant \ $\km^b (f)$ \ (cf.\ 3.2) 
of a \ $\km$-Brunnian link map \ $f$ \ is precisely as
strong as
$$
\widetilde\mu_{\sssize M} (f) \ := \ \{ h_{(g), \gamma} (\km' (f))\} 
\in \ \bigoplus_{((g), \gamma) \in
\Bbb Z^{r -1} \times \Sigma_{r -2}} \ \pi^S_{|p| - (r -1)(m -2) -1}
$$
(see 2.24 as well as \cite{10}, 2.11 and 5.6) on one hand, 
and as the sequence
$$
\quad \left\{ \mum^{(s)} (f) \ \in  \ \ \ \ \bigoplus_{\gamma 
\in \Sigma_{r, s}}
\ \ \ \ \ \pi^S_{|p| - (r -1)(m -2)- 1} \right\}_{s \ge 0}
$$
on the other hand (cf.\ 3.5). 
The transition to the base point free homotopy
invariant \ $\km (f)$ \ amounts to dividing out the translating action
of \ $\Bbb Z^{r -1}$ on the indices \ $((g), \gamma)$ of the factors in
\ $\widetilde\mu_{\sssize M} (f)$ (see \cite{10}, 3.3 and 5.6). 
On the other hand, \; $\mum^{(s)}
(f) (= \mu_{\sssize{\Bbb R^m}} (f^{(s)})$ \ is invariant under
base point free link homotopies as soon as all \ $\mu_{\sssize{\Bbb R^m}}
(f^{(t)}), \ t < s$, are successively defined 
and trivial (or, equivalently, the \lq\lq preceding\rq\rq\ higher \
$\mu$-invariants \ $\mum^{(0)} (f), \dots, \mum^{(s -1)} (f)$ \ vanish;
compare 3.8).
\endremark

\vskip 7 mm
\specialhead  4. \ Connected sums and the extended definition of the invariants
\endspecialhead

It is often possible and useful to define $\kappa$- and $\mu$-invariants for arbitrary (not necessarily $\kappa$-Brunnian) link maps.

Assume that \ $1 \le p_1, \dots, p_r \le m - 3$. Consider a smooth connected target manifold of the form \ $M^m = M' \times \Bbb R$  \ where the base points \ $y^0_1, \dots, y^0_r$ \ lie in \ $M' \times \{ 0\}$. Then there is a well-defined \lq\lq connected sum\rq\rq\  addition on the set \ $BLM_{(p)} (M)$ of base point preserving link homotopy classes of link maps
$$
f \ = \ \amalg f_j \ : \ \amalg S^{p_j} \ \longrightarrow \ M \ \qquad :
\tag 4.1
$$
just deform the summands \ $f^+$ and $f^-$ \ into \ $M' \times [0, \infty)$ \ and \ $M' \times (- \infty, 0]$, resp., and form the obvious sum \ $f^+ + f^-$. If even \ $M = M'' \times \Bbb R^2$ \ this operation makes \ $BLM_{(p)} (M)$ \ into an associative commutative semigroup  with unit (provided \ $\pi_1 (M)$ \ is abelian whenever \ $\min \{p_i\} = 1$). In this case successive addition of sub-link maps of the form
$$
f_{i_1} \ \amalg \dots \amalg \ f_{i_{s -1}} \ \amalg \ f_{i_s} \circ \ \text{reflection} \ \quad, \quad s < r,
$$
(together with \ $r - s$ \ constant component maps) allows us to \lq\lq retract\rq\rq\ (in a $\kappa$-theoretical sense)  arbitrary link maps to \ $\kappa$-Brunnian ones and then to apply our invariants. For details compare \cite{9}, pp. 311--312.

We obtain in particular

\proclaim{Proposition 4.2} \ 
Assume \ $1 \le p_1, \dots, p_r \le m - 3$ \ and \ $1 \le n < m$. Then there are additive invariants
$$ \km' (f) \ \in \ \pip \ \left(S^n \vee \bigvee^{r - 1}_{j = 1} S^{m -1}\right)
$$
and \ $\ \mum^s (f) := H_s (\km' (f)), s = 0, 1, \dots,$ \ canonically defined for \ {\rm every} \ base point preserving link map \ $f \ : \ \amalg S^{p_j} \ \longrightarrow \ M = S^n \times \Bbb R^{m -n}$. We have:
\roster
\item"(i)" \ if $f$ is $\km$-Brunnian, these extended definitions agree with the previous ones (cf.\ {\rm 3.2} and {\rm 3.4});
\item"(ii)" \ these invariants vanish whenever at least one of the components maps \ $f_j$ of $f$ is constant.
\endroster
\endproclaim
The case \ $n = m -1$ \ is taken care of by the isomorphisms \ $BLM_{(p)} (S^{m -1} \times \Bbb R) \ \cong BLM_{(p)} ((S^{m -1} - \{ *\}) \times \Bbb R) \cong BLM_{(p)} (\Bbb R^m)$.
\bigskip

Recall also that in the simply-connected case \ $n \ge 2$ \ base point preserving and base point free link homotopy theory coincide here.
\vskip3mm
Next we extend theorem 3.8.

\proclaim{Theorem 4.3}
\ Under the assumption of proposition $4.2$ the identity
$$
\mug^{(s)} (f) \ := \ H_{s, \gamma} (\km' (f)) \ = \ \murg (f^{(s)})
$$
holds $($at least up to a fixed sign$)$ for {\rm every}\ link map \ $f  :  \amalg S^{p_j} \longrightarrow M  =  S^n \times \Bbb R^{m - n}$ \ $($not necessarily \ $\kappa$-Brunnian; but base point preserving if \ $n = 1)$ \ and for all \ $s \ge 0$ \ and \ $\gamma \in \sum_{s + r - 2}$.
\endproclaim

\demo{Proof} \
Adding suitable link maps in \ $M$ \ with at least one constant component changes neither \ $\mum^{(s)} (f)$ \ nor \ $\mur (f^{(s)})$. Thus we may assume that \ $f$ \ is $\km$-Brunnian. In order to make also \ $f^{(s)}$ \ $\krm$-Brunnian we habe to add successively further link maps in \ $\Bbb R^m$ of the form
$$
\coprod^s_{j = 1} \eta_{z_j} \amalg \ \eta \circ f_1 \amalg \dots \amalg \eta \circ
f_{r -1} \ \amalg \ \eta \circ f_r \circ \text{(reflection)}^i \ ,  \ \ i = 0 \ \text{or}\ 1 \ ,
$$
(cf.\ 3.7), but with at least one component \ $\eta_{z_j}$ replaced by a constant map or, equivalently, the remaining components deformed in such a way that they do not intersect the slice \ $\eta (\{z_j\} \times B^{m -n}) \subset \ \Bbb R^m$ \ (cf.\ 3.6). This changes neither \ $\mum^{(s)} (f) = h \circ (\pinch_s \vee \id) (F_{s+r}^{(s)})$ \ (cf.\ 3.9) nor \ $\mur (f^{(s)})$, and as in the proof of 3.8 above these two terms agree by the projectability theorem 5.2 in 
\cite{9}. \hfill $\blacksquare$
\enddemo

\vskip7mm

\specialhead 5. \ \ Linking coefficients and stable suspensions
\endspecialhead

\proclaim{Theorem 5.1} \
Assume that \ $1 \le p_1, \dots, p_r \le m - 3$ \ and \ $1 \le n \le m - 1$. Let \ $f = \amalg^r f_j : \amalg^r S^{p_j} \longrightarrow M = S^n \times \Bbb R^{m -n}$ \ be a link map (base point preserving if \ $n = 1$) \ such that \ $f_1 \amalg \dots \amalg f_{r -1}$ \ is an embedding (and hence \ $f_r$ \ determines the linking coefficient
$$
\lambda (f) \ \in \ \pi_{p_r} \left( S^n \vee \bigvee^{r -1}_{j = 1} S^{m - {p_j} - 1}\right)  ).
$$
Then the following (systems of) invariants contain an equal amount of information:
\roster
\item"(i)" \ $\km' (f)$ \qquad ;
\item"(ii)" \ $\{ \mum^{(s)} (f) = H_s (\km' (f)\}_{s \ge 0}$ \qquad ;
\item"(iii)" \ $\{ H_s (\lambda (f))\}_{s \ge 0} $ \qquad ;
\item"(iv)" \ the stable suspensions of all the Hilton components of \ $\lambda (f)$ \ corresponding to those basic Whitehead products which involve each meridian sphere \ $S^{m - {p_j} -1}$ \ precisely once, \ $j = 1, \dots, r -1$.
\endroster
\endproclaim

\demo{Proof} \ It is a well-known consequence of the Thom isomorphism and Whitehead theorems that the inclusion \ $\incl$ \ of the wedge
$$
W \ := \ S^n \ \vee \ \bigvee^{r - 1}_{j = 1} \ S_j^{m - {p_j} -1}
\tag 5.2
$$
(formed essentially by a core \ $S^n \times \{*\} \subset S^n \times \Bbb R^{m -n}$ \ and by meridians to \ $\underline f := f_1 \amalg \dots \amalg f_{r -1})$ \ into the link complement of \ $\underline f$ \ induces isomorphisms of homotopy groups up to dimension \ $m - 3$ \ (compare also \cite{10}, 4.6). The Hilton components mentioned above in (iv) constitute the \ $\pi''_* (W)$-part \ $\lambda'' (f)$ \ of $\lambda (f) := \incl_*^{- 1} ([f_r])$ \ (cf.\ 2.9, 2.15, and diagramm 5.3 below).

After a suitable isotopy the images of \ $f_1, \dots, f_{r -1}$ \ intersect \ $S^n \times \Bbb R^{m - n - 1} \times [0, \infty)$ \ in small disjoint half-spheres and we may assume that \ $W$ \ lies in \ $S^n \times \Bbb R^{m -n -1} \times [0, \infty)$. Thus the link homotopy class \ $[f]$ \ is the sum (cf.\ \S\ 4) of \ $[f_1 \amalg \dots, \amalg f_{r -1} \amalg \text{constant}]$ \ with the class \ $e_* (\lambda (f))$ \ consisting of \ $r -1$ \ standard spheres in parallel hyperplanes in a suitable coordinate neighbourhood and of \ $\lambda (f)$ \ (cf.\ section 4 of \cite{10}). Moreover the \lq\lq retraction\rq\rq\ procedure discussed in \S\ 4 has just the effect of replacing \ $\lambda (f)$ \ by its \ $\pi'_* (W)$-part \ $\lambda' (f)$ \ which involves each meridian at least once (compare 2.3). By proposition 4.2 we conclude
$$
\km' (f) \ = \ \km' \circ e_* (\lambda (f)) \ = \ \km' \circ e_* (\lambda' (f)) \ .
$$

For \ $n \ge 2$ \ the proof of our theorem follows from diagram 5.3. Here \ $\Sigma w_{\ell *}, \Sigma \overline w_{\ell *}$
\vskip3mm
$$
\EPSFxsize=12.5cm
\BoxedEPSF{diag53.eps}
$$
\vskip7mm
\flushpar
and \ $H$ \ denote (corresponding) Hilton isomorphisms and Hopf homomorphisms; moreover,
$$
k_s \ := \ s (n -1) + (r - 1) (m -2) - |p| + p_r + 1 \ .
$$
Theorem 4.10 in \cite{10} and results 2.11 and 2.12 above imply that actually
$$
\km' (f) \ = \ \km' \circ e_* (\lambda'' (f))
$$
and that each part of diagram 5.3 commutes up to an automorphism of the target group (which is compatible with the $s$-grading). In particular, for each \ $s \ge 0$ \ the invariants
$$
\mum^{(s)} (f) \qquad \qquad \text{and} \qquad \qquad H_s (\lambda (f)) \ = \ H_s (\lambda'' (f))
$$
are related by an invertible matrix with integer coefficients.

A similar diagram, based on theorem 2.27, yields a proof of our claim also in the case \ $n = 1$. \hfill $\blacksquare$
\enddemo
\bigskip
In view of the last theorem nontriviality, surjectivity, or injectivity results for the stable suspensions \ $E^\infty$ \ have corresponding implications for our invariants. In particular, in the stable dimension range we obtain for the setting of 5.1

\proclaim{Corollary 5.4} \ If \ $p_r/2 \le \sum^{r -1}_{j =1} (m - p_j - 2)$, then $\lambda'' (f)$ \ $($i. e. the \ $\pi'' (W)$-part of the linking coefficient \ $\lambda (f)$, cf.\ {\rm 2.9}\ and {\rm 2.15}$)$ is invariant under link homotopies $($assumed to preserve base points in case \ $n = 1)$ and contains precisely as much information as \ $\km' (f)$.

If \ $n \ge 2$, then for every natural number \ $s$ \ satisfying
$$
\frac{p_r}{2} \ \le \ s (n -1) \ + \ \sum^{r -1}_{j =1} (m - p_j - 2)
$$
the \ $\pi'' (W; s)$-part $($cf. {\rm 2.9}$)$ of \ $\lambda (f)$ \ is a $($base point free$)$ link homotopy invariant of \ $f$ \ and precisely as strong as \ $\mum^{(s)} (f)$; moreover \ $\mum^{(s)}$ \ assumes all values in its target group $($cf. {\rm 3.4}$)$.
\endproclaim 

If \ $f$ \ is \ {\it homotopy Brunnian} \ (i.e.\ every proper sub-link map is nulhomotopic), then under certain dimension conditions the only relevant contribution comes from \ $\lambda'' (f)$ \ (cf.\ 2.28 and compare \cite{10}, 4.6).

\proclaim{Corollary 5.5} \
Assume \ $p_1, \dots, p_r \le m - 3$ \ and
$$
\aligned
|p| \ &\le \ (r - 1)(m - 2) + \frac{p_r}{2} \qquad \qquad \qquad \text{and}\\
|p| \ &\le \ r (m - 2) - p_j \qquad \qquad \qquad, \qquad j = 1, \dots, r - 1.
\endaligned
$$
Then the (base point preserving) link homotopy class of a homotopy Brunnian  link map \ $f : \ \coprod^r_{j = 1} S^{p_j} \longrightarrow M = S^n \times \Bbb R^{m -n}$ \ which embeds \ $\coprod^{r -1}_{j = 1} S^{p_j}$ \ is completely determined by \ $\km' (f)$ \ or, equivalently, by \ $\{ \mum^{(s)} (f)\}_{s \ge 0}$. 

If \ $n \ge 2$, this establishes an isomorphism from the semigroup of such classes onto \ $\bigoplus_{s \ge o}$ \ $(\pi^S_{|p| - s(n - 1) - (r - 1)(m - 2) - 1})^{u_{r, s}}$. 
\endproclaim

Finally we prove our claims concerning the examples 1.5 and 1.6 in the introduction. To begin with, assume \ $p_1 = \dots = p_r = m - 3 = 3$.

Fix at first \ $n \ge 2$ \ and put \ $M = M_n$. Then
$$
\bfmu (f) = \left( \{[f_\ell]\}, \{ (\mum^{(0)} \oplus \mum^{(1)})(f_k \amalg f_\ell)\}, \{ \mum^{(0)} (f_j \amalg f_k \amalg f_\ell)\}\right)
$$
lies in \ $(\pi_3 (S^n))^r \oplus (\Bbb Z_2 \oplus \pi^S_{2 -n})^{r \choose 2} \oplus (\pi_0^S)^{r \choose 3}$; here the first \ $\mum^{(0)}$ \ measures just the \
$\alpha$-invariants (= generalized linking numbers in \ $\pi^S_1 = \Bbb Z_2$) of all 2-component sub-link maps of $f$, when included into \ $\Bbb R^6$ (cf.\ 3.8 and \cite{10}, 2.14).

After a suitable link homotopy we may assume that \ $f$ \ is an embedding (add a small copy of \ $f_j$ \ in a small $6$-ball, $j = 1, \dots, r$, and apply the Whitney trick). As in the proof of theorem 5.1 we obtain therefore
$$
[f] = [f_1 \amalg \dots \amalg f_{r -1} \amalg *] + e_* (\lambda_1 (f)) + e_* (\lambda_2 (f)) + e_* (\lambda_2' (f)) + e_* (\lambda_3 (f))
$$
where \ $f_r$ \ determines the linking coefficient \ $\lambda (f)$ \ in
$$
\pi_3 (S^n \vee \bigvee^{r -1} S^2) \cong \pi_3 (S^n) \oplus \bigoplus^{r -1}_{k = 1} (\pi_3 (S^2) \oplus \pi_3 (S^{n +1})) \oplus \bigoplus_{1 \le j < k < r} \pi_3 (S^3)
$$
and \ $\lambda_1 (f), \dots$ \ denote its Hilton components corresponding to the basic Whitehead products \ $ \iota_0, \{\iota_k\}, \{[\iota_0, \iota_k]\}$ \ and \ $\{[\iota_j, \iota_k]\}$, resp. Their values under \ $e_*$ \ are detected by \ $\bfmu$ \ (cf.\ proposition 2.11 and \cite{10}, theorem 4.10) and yield the contribution for \ $\ell = r$. Note here that, given \ $u \in \pi_3 (S^2) \cong \Bbb Z, \ e_* (\iota_k (u))$ \ lies in an open ball \ $B \ \subset \ S^n \times \Bbb R^{6 - n} \ \subset \ \Bbb R^6$ \ and recall the isomorphism
$$
\alpha \ : \ LM_{3,3} (\Bbb R^6) \ @>{ \ \cong \ }>> \ \pi^S_1 \ = \ \Bbb Z_2
$$
(cf.\ \cite{8}, proposition F, or \cite{4}, theorem 1) which distinguishes \ $e_* (\iota_k (u))$ \ according to the stable suspension (or, equivalently, parity) of \ $u$.

Induction over \ $r$ completes the proof for \ $n \ge 2$.

If \ $n = 1$, approximate \ $f$ \ by a self-transverse immersion and add the nulhomotopic link maps \ $* \amalg \dots \amalg \ \text{reflection} \circ f_j \amalg * \amalg \dots \amalg *, \ j = 1, \dots r$; then we can cancel corresponding double points without complications arising from \ $\pi_1 (M_1)$. The linking coefficient \ $\lambda (f)$ \ of the resulting embedding lies now in the group
$$
\pi_3 ( S^1 \vee \bigvee^{r -1} S^2) \ \cong \ \bigoplus^{r -1}_{k = 1} \left(\bigoplus_{g \in \Bbb Z} \pi_3 (S^2)\right) \oplus \bigoplus_{1 \le j < k < r} \left(\bigoplus_{(g, g') \in \Bbb Z^2} \pi_3 (S^3)\right) \oplus \cdots
$$
with basic Whitehead products \ $\iota_{k, g}$ (cf.\ 2.13 -- 2.14), \ $[\iota_{j, g}, \iota_{k, g'}]$ \ and \ $[\iota_{j, g}, \iota_{j, h}], \ g < h \in \Bbb Z$. The products of the last type can be neglected since \ $e_* ([\iota_{j, g}, \iota_{j, h}])$ \ is trivial by proposition 4.13 in \cite{10}. Also the contributions coming from \ $e_* \circ \iota_{j, g}$ \ factor through the stable suspension homomorphism from \ $\pi_3 (S^2)$ \ to \ $\Bbb Z_2$ \ again and then are detected by $\bfmu$ (cf.\ theorem 2.27 and \cite{10}, theorem 4.10). So are the remaining Hilton parts of \ $\lambda (f)$ \ which yield 3-component sub-link maps of \ $f$. By induction over \ $r$ \ we see that each pair \ $k < \ell$ \ (and each triple \ $j < k < \ell$, resp.) of indices between 1 and $r$ contributes the direct summand
$$\bigoplus_{g \in \Bbb Z} \ \Bbb Z_2 \quad \left(\text{and} \ \ \bigoplus_{(g, g') \in \Bbb Z^2} \ \Bbb Z \quad , \ resp. \right)
$$
to \ $BLM_{3, \dots, 3} (M_1)$. The inclusion \ $M_1 \subset M_2$ \ induces the identification of the meridians \ $\iota_{j, g}, \ g \in \Bbb Z$, with the single meridian \ $\iota_j, \ j = 1, \dots, r$, and hence the summation homomorphisms on \ $\Bbb Z_2^{\infty}$ \ and \ $\Bbb Z^{\infty}$. \hfill $\blacksquare$

The claim in  example 1.6 in the introduction follows from theorem 5.1 above and from proposition 4.13 in \cite{10}.

\subhead Acknowledgements  \endsubhead \ It is a pleasure to thank Uwe Kaiser for stimulating discussions. This work was supported in part within the German--Brazilian Cooperation by IB--BMBF and CNPq.

\vskip1cm

\head{References}
\endhead

\enddocument